\documentclass[12pt,reqno]{amsart}

\usepackage{amssymb}
\usepackage{pdfpages}

\newcommand{\ds}{\displaystyle}     

\DeclareMathAlphabet{\E}{U}{eus}{m}{n}     

\newcommand{\V}{{\mathcal V}}      
\newcommand{\A}{{\mathcal A}}      
\newcommand{\R}{{\mathcal R}}

\newcommand{\PP}{{\mathbb P}}

\newcommand{\N}{{\mathbb N}}

\newcommand{\kk}{{\Bbbk}}

\newcommand{\Q}{{\mathfrak Q}}

\newcommand{\la}{\langle}
\newcommand{\ra}{\rangle}

\newtheorem{thm}{Theorem}[section]
\newtheorem{lemma}[thm]{Lemma}
\newtheorem{cor}[thm]{Corollary}
\newtheorem{prop}[thm]{Proposition}

\theoremstyle{definition}
\newtheorem{defn}[thm]{Definition}

\newtheorem{example}[thm]{Example}

\newtheorem{ques}[thm]{Question}

\newtheorem{Pf}{Proof$\!\!$}         
         
\newenvironment{pf}{\begin{Pf}}{\qed\end{Pf}}

\renewcommand{\qedsymbol}{{\hspace*{\fill}\rule{2.7mm}{2.7mm}}}

\DeclareMathSymbol{\twoheadrightarrow}  {\mathrel}{AMSa}{"10}

\newcounter{letter}
\renewcommand{\theletter}{\rom{(}\alph{letter}\rom{)}}

\newcounter{rnum}
\renewcommand{\thernum}{\rom{(}\roman{rnum}\rom{)}}

\begin{document}

\title[On the Notion of Complete Intersection]%
{On the Notion of Complete Intersection outside\\[2mm] 
the Setting of Skew Polynomial Rings}


\maketitle

\vspace*{0.1in}

\baselineskip15pt

\begin{center}
\renewcommand{\thefootnote}{\fnsymbol{footnote}}
{\sc Michaela Vancliff}\footnote{M.~Vancliff was
supported in part by NSF grants DMS-0900239 and DMS-1302050.}\\
Department of Mathematics, P.O.~Box 19408,\\
University of Texas at Arlington,\\
Arlington, TX \ 76019-0408\\
{\sf vancliff@uta.edu \quad uta.edu/math/vancliff}
\end{center}

\setcounter{page}{1}
\thispagestyle{empty}

\bigskip
\bigskip

\begin{abstract}
\baselineskip13pt
In recent work of T.\ Cassidy and the author, a notion of complete 
intersection was defined for (non-commutative) regular skew polynomial rings, 
defining it using both algebraic and geometric tools, where the commutative 
definition is a special case.
In this article, we extend the definition to a larger class of algebras
that contains regular graded skew Clifford algebras, the
coordinate ring of quantum matrices and homogenizations of universal
enveloping algebras.
Regular algebras are often considered to be non-commutative analogues of
polynomial rings, so the results herein support that viewpoint. 
\end{abstract}

\baselineskip18pt


\bigskip

\section*{Introduction}

In \cite[Definition~13]{CVc}, a notion of complete intersection was 
given for (Artin-Schelter) regular skew polynomial rings, defining it 
using both algebraic and geometric tools. That definition is modeled 
on the commutative definition and modified for the context of skew 
polynomial rings.  In this article, we prove that the definition extends
to a larger family of algebras; a family that contains the polynomial ring, 
regular skew polynomial rings, regular graded Clifford algebras, 
regular graded skew Clifford algebras, the coordinate ring of quantum
matrices, and homogenizations of universal enveloping algebras of 
finite-dimensional Lie algebras. 

As in \cite[Definition~13]{CVc}, the definition of complete intersection 
for polynomial rings is a special case of our definition. However, it is
worth noting that other notions of complete intersection abound in the 
literature, with most emphasizing a homological approach, such as the recent 
work in~\cite{KKZ}.

The main focus in this article are graded skew Clifford 
algebras (GSCAs). Such algebras were first defined in \cite{CV}, and the 
family of regular GSCAs contains the family of regular skew polynomial 
rings and regular graded Clifford algebras. In fact, GSCAs may be considered a 
quantized analogue of graded Clifford algebras, and they are of interest in the 
classification of quadratic regular algebras, in part because all regular
algebras of global dimension at most two are GSCAs and in part since most 
quadratic regular algebras of global dimension three may be generated using 
GSCAs (\cite{CV,NVZ}). 
Hence, GSCAs seem a natural place to start in an attempt to extend the 
definition of complete intersection to a setting outside that of 
skew polynomial rings.
We also consider the notion of complete intersection for algebras
that are not GSCAs.

In Section~1, we give the definition of a GSCA and
list the main results of \cite{CV,CVc} that will be useful herein.
In \S\ref{conds}, we list the conditions of interest to us in 
our proposed definition (Definition~\ref{compintn2}) of complete intersection
to follow in Section~3; 
these conditions entail three algebraic tools and one geometric tool.
In Section~2, we consider the definition in the context of certain algebras 
that are related to skew polynomial rings and prove in Theorem~\ref{thm1} 
that most of the definition is applicable there. We also apply 
Theorem~\ref{thm1} to some regular algebras that are not GSCAs.

Our main results are given in Section~3, in Theorems~\ref{thm2} and
\ref{thm3},
where we prove that our proposed definition is applicable to regular GSCAs
and to a larger family of algebras that contains the coordinate ring of
quantum matrices.
Regular algebras are often considered to be non-commutative analogues of
polynomial rings, so these results support that viewpoint. It is therefore of
interest to see if our definition of complete intersection (or a modification 
thereof) is applicable to a more general family of algebras than
those considered in Theorems~\ref{thm2} and~\ref{thm3}.

\bigskip
\bigskip


\section{Graded Skew Clifford Algebras} 
In \S\ref{gsca}, we define graded skew Clifford algebras and list 
the results of \cite{CV,CVc} that will be of interest in the subsequent
sections.  We list,
in \S\ref{conds}, some conditions that will be useful in our 
proposed definition (Definition~\ref{compintn2})
of complete intersection in Section~3.

Throughout the article, $\kk$ denotes an algebraically closed field with 
char$(\kk) \neq 2$, $M(n,\ \kk)$ denotes the space of $n \times n$ matrices 
with entries in $\kk$, and, for any vector space $V$, we let $V^\times = 
V\setminus \{ 0 \}$. 

\medskip

\subsection{The Algebras and Known Results}\label{gsca}\quad

The following notion of symmetry for square scalar matrices generalizes the
standard notions of symmetry and skew-symmetry for such matrices.

\begin{defn}\cite{CV}\label{mu-symm}
For $\{i, j \}\subset \{ 1, \ldots , n\}$, let $\mu_{ij} \in \kk$ 
satisfy $\mu_{ij} \mu_{ji} = 1$ for all $i \neq j$, and 
write $\mu = (\mu_{ij}) \in M(n,\ \kk)$.  
A matrix $M \in M(n,\ \kk)$ is called $\mu$-{\em symmetric} if
$M_{ij} = \mu_{ij}M_{ji}$ for all $i, j = 1, \ldots , n$. 
We denote the space of such matrices by $M^\mu(n,\ \kk)$.
\end{defn}
\begin{defn}\cite{CV}\label{GSCA}
With $\mu$ as above, suppose $\mu_{ii} = 1$ for all $i$, and fix 
$\mu$-symmetric matrices 
$M_1, \ldots , M_n \in M^\mu(n,\ \kk)$.  A {\em graded skew Clifford algebra} 
(GSCA)
associated to $\mu$ and $M_1$, $\ldots ,$ $M_n$ is a graded $\kk$-algebra
on degree-one generators $x_1, \ldots , x_n$ and on degree-two generators
$y_1, \ldots , y_n$ with defining relations given by:
\begin{enumerate}
\item[(a)] (degree-2 relations)
$\ds x_i x_j + \mu_{ij} x_j x_i = \sum_{k=1}^n (M_k)_{ij} y_k$
           for all $i, j = 1, \ldots , n$, and
\item[(b)] degree-3 and degree-4 relations that guarantee the existence of a 
           normalizing sequence $\{ r_1, \ldots , r_n\}$
           of homogeneous elements of degree two that span 
	   $\kk y_1 + \cdots + \kk y_n$. 
\end{enumerate}
\end{defn}
\noindent 
Clearly, if $\mu_{ij} = 1$ for all $i, j$, and if 
$y_k$ is central in the algebra for all~$k$, then the algebra is a graded 
Clifford algebra. 
On the other hand, taking $M_1, \ldots , M_n$ to be diagonal and linearly 
independent shows that regular skew polynomial rings on finitely-many 
generators are GSCAs.

Although GSCAs need not be quadratic nor regular
in general, a simple geometric criterion was established in 
\cite[Theorem~4.2]{CV} for determining when such an algebra is quadratic and 
regular. In order to state that result, we first associate a (non-commutative)
quadric system to a GSCA. 

\begin{defn}\cite{CV,CVc}\label{quadric}
Fix $\mu$ as in Definition~\ref{GSCA}, and let $S$ denote the regular skew 
polynomial ring 
\[
S= 
\frac{
T}
{\langle  z_j z_i - \mu_{ij} z_i z_j \ : \ i, j = 1, \ldots , n\ \rangle },
\\[3mm]
\] 
where $T$
is the tensor (free) algebra on generators $z_1, \ldots , z_n$.  
A (noncommutative) quadratic form is defined to be any element of $S_2$.
By \cite{CV}, we have $M^\mu(n,\ \kk) \cong S_2$, as vector spaces, via 
$M \mapsto z^T M z \in S$, where  $z = (z_1, \ldots , z_n)^T$. This map
mirrors the isomorphism between symmetric matrices and commutative quadratic 
forms.
The span of quadratic forms $Q_1, \ldots , Q_m\in S_2$ will be called the 
{\em quadric system} associated to $Q_1, \ldots , Q_m$.  If a quadric system 
is given by a normalizing sequence in $S$, then it is called a {\em normalizing
quadric system}.  
\end{defn}

The following uses the notion of a 1-{\em critical} module (with respect
to GK-dimension) from~\cite{Lev}, which is a module~$M$ such that GKdim$(M) = 
1$ and GKdim$(M/N)< 1$ for all nonzero submodules~$N$ of~$M$; if, further, 
$M$~is graded, it suffices to consider only nonzero graded submodules of~$M$.

\begin{defn}\cite{CV,CVc}\label{BP}
We define a {\em right base point} of a quadric system 
$\Q \subset S_2$ to be any right base-point module over $S/\la \Q \ra$; that 
is, to be any 1-critical graded right module over $S/\la \Q \ra$ that is 
generated by its homogeneous degree-zero elements and which has Hilbert series 
$H(t) = c/(1-t)$, for some $c \in \N$.  We say a quadric system is 
{\em right base-point free} if it has no right base points. Analogous
definitions hold for left modules, etc.
\end{defn}

One may associate some geometry to the definitions of base point and
base-point free as follows. By identifying $\PP(S_1^*)$ with $\PP^{n-1}$, 
let $Z$ denote the zero locus, in $\PP^{n-1} \times \PP^{n-1}$, of the defining
relations of $S$.
For any $q \in S_2$, we call the zero locus of~$q$ in~$Z$ the {\em quadric 
associated to}~$q$, and denote it by $\V_{_S}(q)$; in other words, 
$\V_{_S}(q) = \V(\hat q) \cap Z$, where $\hat q$ is any lift of $q$ to 
$T_2$, and $\V(\hat q)$ is the zero locus of $\hat q$ in 
$\PP^{n-1} \times \PP^{n-1}$.  By \cite{ATV1}, if $\Q$ is a normalizing quadric
system, then the isomorphism classes of right (respectively, left) base-point 
modules $M=\bigoplus_{i=0}^\infty M_i$ over $S/\la \Q \ra$ such that 
$\dim_\kk(M_0) = 1$ are
parametrized by the scheme $\bigcap_{q\in\Q} \V_{_S}(q)$ (since the
projection of $Z$ to the left $\PP^{n-1}$ equals the projection of $Z$ to the 
right $\PP^{n-1}$ and by applying \cite[Lemma~3.2]{CV}).

If $S$ is commutative, then $\dim_\kk(M_0) =1$ for all base-point 
modules, so, in the commutative setting, the previous paragraph implies that 
the above definition of ``base point'' agrees with the commutative definition 
thereof.

By \cite[Corollary~11]{CVc}, a normalizing quadric system $\Q$ is right 
base-point 
free if and only if $\dim_\kk(S/\la \Q \ra) < \infty$. Hence, such a quadric 
system is right base-point free if and only if it is left base-point free.
In particular, the adjectives ``right'' and ``left'' may be dropped when 
referring to a {\em normalizing} quadric system being base-point free.

\begin{thm}\label{cvthm}\cite{CV,CVc}
A graded skew Clifford algebra $A = A(\mu, M_1, \ldots , M_n)$ as in
Definition~\ref{GSCA} is a quadratic, Auslander-regular algebra of global 
dimension $n$ that satisfies the Cohen-Macaulay property with Hilbert series 
$1/(1-t)^n$ if and only if the quadric system associated to 
$M_1 , \ldots , M_n$ is normalizing and base-point free; in
this case, $A$ is a noetherian regular domain and is unique up to 
isomorphism.
\end{thm}

The ``base-point free'' property in the previous result is referring to the 
associated quadric system in the skew polynomial ring $S$ yielding a complete
intersection in the sense of Definition~\ref{compintn} below. 

\newpage

\subsection{The Conditions}\label{conds} \quad

Let $A = \bigoplus_{i=0}^\infty A_i$ denote an $\N$-graded,
finitely generated $\kk$-algebra. The following definition is motivated
by Definition~\ref{BP}.

\begin{defn}\label{D}
Suppose further that $A$ is connected and generated by $A_1$.
We define a right (respectively, left) base-point module 
over $A$ to be any 1-critical (with respect to GK-dimension) graded 
right (respectively, left) $A$-module that is generated by its 
homogeneous degree-zero elements and which has Hilbert series 
$H(t) = c/(1-t)$ for some $c \in \N$.
\end{defn}

Suppose that GKdim$(A) = n \in \N$ and that $A$ has a normalizing sequence 
$F=\{ f_1 , \ldots , f_t\} \subset A\setminus\kk^\times$ consisting of 
homogeneous elements.
In order to address the notion of complete intersection in the context of 
$A$ and $F$, we will be interested in four (possibly equivalent) conditions on 
$A$ and $F$ as follows:
\begin{enumerate}
\item[I.] $F$ is a regular sequence in $A$;
\item[II.] $\ds \dim_\kk(A/\la F\ra) < \infty$;
\item[III.] (if $t \leq n$) for each $k = 1, \ldots , t$,
           GKdim$\left( A/\la f_1 , \ldots , f_k\ra \right) = n-k$;
\item[IV.] there are no right base-point modules over $A/\la F \ra$ (if 
           $A/\la F \ra$ is connected and generated by its degree-1
	   elements). 
\end{enumerate}
By symmetry of conditions I-III, if condition IV is equivalent to any of
the other conditions, then it will hold if and only if there are no 
left base-point modules over $A/\la F \ra$. Hence, unless stated otherwise, 
we will focus on right base-point modules.

If $A$ is the algebra $S$ from Definition~\ref{quadric} and 
if $\deg (f_k) = 2$ for all $k$, then condition~IV is equivalent to the 
normalizing quadric system $\{ f_1, \ldots , f_n\}$ being base-point free.

\bigskip

The definition of complete intersection in \cite{CVc} that we wish to 
extend to regular GSCAs is based on 
\cite[Corollary~11]{CVc}, which proves that conditions I-IV are
equivalent in the case that $A$~is the algebra $S$ from
Definition~\ref{quadric}.
This motivates the following definition in \cite{CVc}.

\begin{defn}\cite{CVc}\label{compintn}
Let $A$ denote the skew polynomial ring $S$ from Definition~\ref{quadric}.
If $F$ is a normalizing sequence in $A$ of $n$ homogeneous elements of positive
degree, then we call $\ds A/\la F \ra$ a {\em complete intersection} if the 
equivalent conditions I-IV hold.
\end{defn}
\noindent
In Definition~\ref{compintn2}, we extend Definition~\ref{compintn} to regular 
GSCAs and a larger family of algebras.

\medskip

Henceforth, we use the term ``regular skew polynomial ring'' to refer to any
algebra having defining relations like those defining $S$ in 
Definition~\ref{quadric}.

\bigskip


\section{Algebras Related to Regular Skew Polynomial Rings}

It is a natural problem to seek the algebras $A$ for which two or more of 
conditions~I-IV from \S\ref{conds} are equivalent.  
In this section, we address this issue for certain types of algebras related
to regular skew polynomial rings.

The technical hypotheses of the next result are motivated by results in 
\cite{Lev}.  As in \cite{Lev}, we assume that the definition of 
Auslander-Gorenstein includes the requirement that the algebra be noetherian.

\begin{thm}\label{thm1}
Let $A$ denote an $\N$-graded Auslander-Gorenstein $\kk$-algebra of finite 
injective dimension that satisfies the Cohen-Macaulay property with
GKdim$(A) = n \in \N\setminus \{ 0\}$.  Suppose $F\subset A$ is a normalizing 
sequence consisting of $n$~homogeneous elements of positive degree.
\begin{enumerate}
\item[(a)] \cite[Theorem~2.5]{CV}
	   Conditions I-III are equivalent.
\item[(b)] If $A/\la F \ra \cong S/\la G\ra$, where $S$ is a regular skew
           polynomial ring on $m$ generators and $G$ is a normalizing sequence 
	   in $S\setminus\kk^\times$ of $m$ homogeneous elements, 
	   then conditions~I-IV are equivalent for $A$ and $F$.
\end{enumerate}
\end{thm} 
\begin{pf} 
(b) \ Suppose $A/\la F \ra \cong S/\la G\ra$, where $S$ and $G$ are as given 
in the statement. 
Condition~IV holds if and only if $S/\la G \ra$ has no right base-point 
modules.  By \cite[Corollary~11]{CVc}, this situation happens if and only if 
$\dim_\kk(S/\la G \ra) < \infty$. It follows that condition~IV is equivalent to 
condition~II.  Hence, the result follows by (a).
\end{pf}

Since many known Auslander-regular (respectively, AS-regular) algebras satisfy 
the hypotheses of Theorem~\ref{thm1}, we have that Theorem~\ref{thm1} may be 
applied to many such algebras. 
In Example~\ref{qumat}, we present an application of Theorem~\ref{thm1},
where the algebra is not a GSCA.

\begin{example}\label{qumat}
Let $A$ denote the coordinate ring of quantum $2 \times 2$ matrices; that is, 
$A$ is the quadratic $\kk$-algebra on generators $a, b, c, d$ with defining 
relations
\[
\begin{array}{ccc}
ab = q ba, &\qquad & bd = q db,\hspace*{2ex}\\[1mm]
ac = q ca, &      & cd = q dc,\hspace*{2ex}\\[1mm]
bc = cb,~  &      & ad - da = (q - q^{-1}) bc,
\end{array}
\]
where $q \in \kk^\times$ and $q^2 \neq 1$ (\cite{FRT}). By \cite{V1}, this 
algebra is an Auslander-regular algebra of global dimension four with Hilbert
series $1/(1-t)^4$ and it satisfies the hypotheses of Theorem~\ref{thm1},  
but it is not a GSCA. Let
$F = \{ bc,\ b^2,\ c^2,\ ad \}$, which is a normalizing nonregular sequence 
in~$A$, so condition~I fails.  Since $A/\la F \ra \cong S/\la G \ra$, 
where $S$ is the $\kk$-algebra on generators $a, b, c, d$ with defining
relations
\[
\begin{array}{ccccc}
ab = q ba, &\qquad & bd = q db,&\qquad & bc = cb,\\[1mm]
ac = q ca, &       & cd = q dc,&       & ad = da,
\end{array}
\]
and $G = \{ bc,\ b^2,\ c^2,\ ad \} \subset S$, conditions~I-IV are equivalent 
for $A$ and $F$, by Theorem~\ref{thm1}. 
\end{example}

In (a)-(c) of the next result, it should be noted that it is {\em not}
assumed that 
the quadric system associated to the GSCA $A$ is linearly
independent, nor that $A$ is quadratic.

\begin{prop}\label{prop1}
Let $A$ denote a graded skew Clifford algebra and let $F$ denote the
normalizing sequence of length $n$ given in Definition~\ref{GSCA}(b).
\begin{enumerate}
\item[(a)] The algebra $A$ is noetherian.
\item[(b)] Conditions II and IV hold for $A$ and $F$.
\item[(c)] If $F$ is regular in $A$, then $A$ satisfies the hypotheses of 
           Theorem~\ref{thm1} and has Hilbert series $H(t) = 1/(1-t)^n$, and 
	   conditions~I-IV hold for $A$ and $F$.
\item[(d)] If $F$ is regular in $A$ and if $A$ is quadratic, then $A$ is
           Auslander-regular of global dimension $n$ and the quadric
	   system associated to $A$ is normalizing and base-point free.
\end{enumerate}
\end{prop}
\begin{pf}
In the notation of Definition~\ref{GSCA}, we have 
\[
\frac{A}{\la F \ra} \cong
\frac{A}{ \sum_{k=1}^n y_k A} \cong
\frac{\kk\la x_1 , \ldots , x_n\ra}{
\la x_i x_j + \mu_{ij} x_j x_i : 1 \leq i,\ j \leq n \ra 
} \cong
\frac{B}{\la G \ra}, 
\]
where $B$ is a regular skew polynomial ring on generators $x_1, \ldots , x_n
\in B_1$ and $G = \{ x_1^2, \ldots , x_n^2\}$.

(a) Since $B/\la G\ra$ is noetherian, it follows that $A$ is noetherian by 
    \cite[Lemma~8.2]{ATV1}.

(b) The algebra $B/\la G\ra$ has finite dimension and so has no base-point
    modules.  Hence, conditions~II and IV hold for $A$ and $F$.

(c) The sequence $G$ is normalizing in $B$ and is base-point free. By
    \cite[Corollary~11]{CVc}, $G$ is a regular sequence in $B$, and, by
    \cite[Theorem~5.10]{Lev}, $B/\la G \ra$ is Auslander-Gorenstein of
    finite injective dimension and satisfies the Cohen-Macaulay property.  
    Thus, if $F$ is regular in $A$, then, by \cite[Theorem~5.10]{Lev}, $A$
    satisfies the hypotheses of Theorem~\ref{thm1}, so conditions~I-IV hold 
    for $A$ and $F$. Since the Hilbert series of $B/\la G\ra$ is $h(t) =
    (1+t)^n$ and since $F$ is regular and normalizing in $A$, the Hilbert 
    series of $A$ is $H(t) = 1/(1-t)^n$. 

(d) By \cite[Theorem~1.5]{S.Tingey}, $B/\la G\ra$ is Koszul, since $B$ is 
    Koszul and $G$ is regular and normalizing. Thus, $A/\la F\ra$ is 
    Koszul. Since $A$ is quadratic and $F$ is regular and normalizing, $A$ is
    Koszul by \cite[Theorem~1.5]{S.Tingey}. By (c), the Koszul dual of $A$ 
    has Hilbert series $h(t) = (1+t)^n$ and $A$ is Auslander-Gorenstein. It 
    follows that the right trivial 
    module~$\kk_A$ has projective dimension~$n$, so $A$ is Auslander-regular 
    of global dimension~$n$. The result follows by Theorem~1.4.
\end{pf}

The following example shows that $A$ being quadratic does not follow
from $F$ being regular in Proposition~\ref{prop1}(c)(d).

\begin{example}
Let $A$ denote the graded Clifford algebra given by the symmetric matrices
\[
M_1 = \begin{bmatrix} 2 & -1 \\ -1 & 0  \end{bmatrix},
\quad 
M_2 = \begin{bmatrix} 0 & -1 \\ -1 & 2  \end{bmatrix}.
\]
Hence, $A$ is a factor algebra of the $\kk$-algebra $A'$ on degree-1 
generators $x,\ y$ with defining relations
\[
xy + yx = -x^2 - y^2, \qquad x^2 y = y x^2.
\]
Since $x^2$ and $y^2$ are central
elements in $A'$, we have that $A = A'$.  Applying Bergman's Diamond
Lemma, the ambiguities are resolvable, so $\dim_\kk(A_m) = m+1$, for all $m
\in \N$, thus $A$ has Hilbert series $H(t) = 1/(1-t)^2$. 
Since $x^2$ and $y^2$ are central and $A/\la x^2,\ y^2\ra$ has Hilbert series 
$h(t) = (1+t)^2$, it follows that $\{x^2,\ y^2\}$ is a regular normalizing 
sequence in~$A$.
Thus, by Proposition~\ref{prop1}(c), conditions I-IV hold for $A$ and
the sequence $F=\{x^2,\ y^2\}$. Although $A$ satisfies the hypotheses of 
Theorem~\ref{thm1}, one should note that $A$ is not quadratic and not
regular (since $(x+y)^2 = 0$ in $A$).
The author thanks James Zhang of the University of Washington for providing 
this example.
\end{example}

\bigskip


\section{Complete Intersections}

Continuing our search for algebras $A$ for which two or more of 
conditions~I-IV from \S\ref{conds} are equivalent, we will prove
in this section, in Theorem~\ref{thm2}, that if $A$ is a regular GSCA of global
dimension~$n$ (in the sense of Theorem~\ref{cvthm})
and if $F$ is any normalizing sequence in~$A$ as in 
\S\ref{conds} of length $n$, 
then conditions~I-IV are equivalent.   We also prove a
similar result, in Theorem~\ref{thm3}, for a family of algebras that contains 
regular GSCAs. 

In order to motivate our approach, we note that the algebra $A$ in
Example~\ref{qumat}
is regular but not a GSCA. However, that algebra contains a normalizing
sequence ${\mathcal F} = \{b, c, d \}$ such that GKdim$(A/\la \mathcal
F\ra ) = 1$.  We will use such a property in Proposition~\ref{P}.

\begin{lemma}\label{L}
Let $B=\bigoplus_{i= 0}^\infty B_i$ denote an $\N$-graded 
$\kk$-algebra such that $\dim_\kk(B_i) < \infty$ for all~$i$. If 
$\dim_\kk(B/bB) < \infty$,  for some $b \in B_d$, where $d > 0$, 
then GKdim$(B) \leq 1$.
\end{lemma}
\begin{pf}
Suppose $\dim_\kk(B/bB) < \infty$, for some $b \in B_d$, where $d > 0$. 
Thus, $B_{k+d} = (bB)_{k+d} = b B_k$ for all $k \gg 0$.  Hence, 
$\dim_\kk(B_{k+d}) \leq \dim_\kk(B_k)$ for all $k \gg 0$.  It follows that 
$B$~is finitely generated and that the Hilbert series of~$B$ is eventually 
periodic, and so GKdim$(B) \leq 1$.
\end{pf}

For the rest of this section, we write $\widehat{\V_A(I)}$  to denote the set 
of isomorphism classes of (right) base-point modules over $A/I$ where
$A= \bigoplus_{i=0}^\infty A_i$ is a connected, $\N$-graded, finitely 
generated $\kk$-algebra generated by $A_1$ and $I$ is any
homogeneous ideal of $A$.

\begin{prop}\label{P}
Let $A= \bigoplus_{i=0}^\infty A_i$ denote a noetherian, connected, 
$\N$-graded $\kk$-algebra that is generated by $A_1$. Suppose that 
there exists a normalizing sequence $\{ y_1, \ldots , y_\nu\} \subset A
\setminus \kk$ consisting of homogeneous elements such that 
GKdim$(A/\la y_1, \ldots , y_\nu \ra ) = 1$. Let 
$F \subset A\setminus \kk^\times$ denote any finite normalizing sequence 
consisting of homogeneous elements. 
\begin{enumerate}
\item[(a)] \  
Conditions II and IV are equivalent for $A$ and $F$; that is,
$\dim_\kk(A/\la F\ra) < \infty$ if and only if \/ 
$\widehat{\V_A(\la F \ra)}$ is empty.
\item[(b)] \  
There is a one-to-one correspondence between the set of isomorphism 
classes of right base-point modules over $A/\la F \ra$ with Hilbert series
$H(t) = c/(1-t)$, where $c \in \N$, and the set of isomorphism classes of 
left base-point modules over $A/\la F \ra$ with Hilbert series 
$H(t) = c/(1-t)$. 
\end{enumerate}
\end{prop}
\begin{pf}
(a) \  
Let $I = \la F \ra$ and suppose $M = \bigoplus_{i=0}^\infty M_i \in 
\widehat{\V_A(I)}$. In particular, $M = M_0 (A/I)$ and $\dim_\kk(M_0) < \infty$,
while $\dim_\kk(M) = \infty$; thus, $\dim_\kk(A/I) = \infty$.

Conversely, suppose $\dim_\kk(A/I) = \infty$.  Setting $y_0 = 0 \in A$,
the hypotheses and Lemma~\ref{L} imply that there exists $m \in \{0,
\ldots, \nu\}$ such that GKdim$(A/(I + \sum_{i=0}^m y_i A)) = 1$.
Let $k$ denote the smallest such $m$, and let 
$\A = A/(I + \sum_{i=0}^k y_i A)$, which is a $\kk$-algebra since 
$\{ y_1, \ldots , y_\nu\}$ is a normalizing sequence in~$A$. By construction, 
GKdim$(\A) = 1$. 
We may now apply \cite[Example~5.5]{AZ} to $\A$ as follows. 

By \cite{SSW}, $\A$ is P.I.\ and finitely generated over its
center and contains a regular homogeneous central element~$z$ of degree~$d$ 
for some positive integer~$d$.  Let $B = \A[z^{-1}]$ and $R = B_0$.
The algebras $\A$ and $B$ are locally finite, so $\dim_\kk(R) < \infty$ and 
so $R$ has a finite-dimensional simple module $N$. Moreover, $B\cong
R[x, x^{-1} ; \sigma]$ for some $x \in B_1$ and $\sigma \in$
Aut$(R)$, by \cite[Corollary I.3.26]{NV}. By \cite[Example~5.4]{AZ},
proj$(\A) \cong $ spec$(R)$, where proj$(\A)$ denotes the category of
finitely generated graded $\A$-modules modulo the subcategory of finitely
generated graded torsion $\A$-modules. In particular, by \cite[\S7]{ATV2},
let $\hat N = \bigoplus_{i=0}^\infty (N \otimes_R B)_i = 
\bigoplus_{i=0}^\infty (N \otimes_R B_i)$, which is a finitely
generated graded $\A$-module, generated by $\hat N_0 = N\otimes\kk$. 
In fact, as vector spaces, $\hat N_i
= N \otimes \kk x^i$. Hence, in order to prove that $\hat N$ is a base-point
module, it suffices to prove it is 1-critical; 
to this end, we follow some of the argument in \cite{CVc}.

By definition, $\hat N_i$ is an $R$-module.
Let $e \in \hat N_i^\times$, so $e = v \otimes x^i$ for some $v \in N^\times$.
Given that $x$ is normal and regular in $B$, we have 
$eR = (v \otimes x^i)R = v \otimes R x^i = vR \otimes \kk x^i = N \otimes
\kk x^i$, since $N$ is a simple $R$-module. Hence, $\hat N_i$ is a simple
$R$-module for all $i$.

Let $M = \bigoplus_{i=r}^\infty M_i$ denote a nonzero graded $\A$-submodule
of $\hat N$, where $M_r \neq \{ 0\}$. Since $z$ acts faithfully on $\hat N$, 
$z$ acts faithfully on $M$, so, for each $i \in \{r, \ldots , r+d-1\}$, we
have $\{\dim_\kk (M_{i+td})\}_{t \ge 0}$ is a
nondecreasing sequence. However, $\dim_\kk(\hat N_j) = \dim_\kk(N)$ for all $j$,
so it follows that $\dim_\kk(M_j) = \dim_\kk(M_{j+d})$, for all $j\gg 0$.
Hence, $M_j$ is a nonzero $R$-submodule of the simple $R$-module $\hat N_j$
for all $j\gg 0$; whence, $M_j = \hat N_j$ for all $j\gg 0$. 

It follows that $\hat N$ is a right base-point module over $\A$, and hence a 
right base-point module over $A/I$, so $\hat N \in \widehat{\V_A(I)}$.

(b) \   
Let $M \in \widehat{\V_A(I)}$, where $\dim_\kk(M_0) = c \in \N$, and 
let $\A = A/\text{Ann}(M)$.
If $y_i \in$ Ann$(M)$ for all $i$, then GKdim$(\A) = 1$, so, 
using the proof in~(a), it follows that, in this case, there exists a
homogeneous regular central element $y \in \A$.  On the other hand, if 
$y_i \notin$ Ann$(M)$ for some $i$, let $k$ denote the least such $i$ and
let $y$ denote the image of $y_k$ in $\A$.  In either case, by
\cite[Lemma~5]{CVc}, $y$ is a regular normal homogeneous element of $\A$. 
Let $d = \deg(y)$, $B = \A[y^{-1}]$ and $R = B_0$.

Since the action of $y$ on $M$ is a bijective linear map from $M_i \to
M_{i+d}$ for all $i$, each $M_i$ may be viewed as an $R$-module. Let $N
\subset M_j$ denote a nonzero $R$-submodule of $M_j$, for some~$j$.
Since $M$ is 1-critical, $N\A_i=(N\A)_{i+j} = M_{i+j}$ for all $i \gg 0$.
However, using the bijective action of $y$, we have 
$M_j = M_{td+j}y^{-t} = N\A_{td}y^{-t} \subset N$ for all $t \gg 0$, and so 
$N = M_j$. Thus, $M_j$ is a simple $R$-module for all $j$.

Let $\R = R/\text{Ann}(M_0)$. Since $\dim(\R) < \infty$, and $M_0$ is a
simple $R$-module, $\R$ is a prime
Artinian ring, and so is a simple ring, and, by the Artin-Wedderburn
Theorem (since $\kk$ is algebraically closed), we have $\R \cong M(c,\ \kk)$. 
As there exists a one-to-one
correspondence between simple right modules over $M(c,\ \kk)$ and simple 
left modules over $M(c,\ \kk)$ (by taking the transpose of the matrices
and vectors), there exists a simple left module $N$ over
$M(c,\ \kk)$ that we may associate to $M_0$.  Let $\hat N =
\bigoplus_{i=0}^\infty (B \otimes_R N)_i$.

However, $M \cong
\bigoplus_{i=0}^\infty (M_0 \otimes_R B)_i$ and has Hilbert series $H_M(t) = 
c/(1-t)$, so it follows that, for each $i$, we have 
$B_i = R \otimes_\kk V_i$ for some 1-dimensional subspace $V_i$ of $B_i$.
Hence, for each $i$, as vector spaces we have $B_i = W_i \otimes_\kk R$
for some 1-dimensional subspace $W_i$ of $B_i$. It follows that 
$\hat N = \bigoplus_{i=0}^\infty (W_i \otimes_\kk N)$, and so $\hat N$ has
Hilbert series $H(t) = c/(1-t)$. Using an argument similar to that in (a),
it follows that $\hat N$ is 1-critical, and so is a left base-point module
over~$\A$. Applying this discussion to $\hat N$ yields $M$, so (b) follows.
\end{pf}

\begin{cor}\label{C}
Let $A$ denote a GSCA as in Definition~\ref{GSCA} and assume that the 
$\mu$-symmetric matrices $M_1, \ldots , M_n$ are linearly independent.
If $F \subset A\setminus \kk^\times$ is any finite normalizing sequence 
consisting of homogeneous elements, then both (a) and (b)
of Proposition~\ref{P} hold for~$A$ and~$F$.
\end{cor}
\begin{pf}
By \cite[Lemma~1.13]{CV}, $A$ is generated by degree-1 elements if and only
if the defining $\mu$-symmetric matrices are linearly independent.
Moreover, $A$ has the property that 
$\dim(A/\la r_1, \ldots , r_n\ra ) < \infty$,
where $\{ r_1, \ldots , r_n\}$ is the normalizing sequence in 
Definition~\ref{GSCA}(b).  Hence, by Proposition~\ref{prop1}(a) and
Lemma~\ref{L}, $A$ satisfies the hypotheses of Proposition~\ref{P},
so the result follows.
\end{pf}

\begin{thm}\label{thm2}
Let $A$ denote a regular graded skew Clifford algebra of global dimension~$n$
in the sense of Theorem~\ref{cvthm}.
If $F\subset A\setminus \kk^\times$ is any normalizing sequence consisting
of $n$~homogeneous elements, then conditions~I-IV are equivalent.
\end{thm}
\begin{pf}
The hypotheses imply that the normalizing sequence from 
Definition~\ref{GSCA}(b) is regular, so that Proposition~\ref{prop1}(c) 
applies to $A$.  Hence, conditions I-III are equivalent by 
Theorem~\ref{thm1}(a). Conditions II and IV are equivalent by
Corollary~\ref{C}.
\end{pf}

More generally, we have the following result by combining Theorem~\ref{thm1}(a) 
with Proposition~\ref{P}(a).

\begin{thm}\label{thm3}
Let $A= \bigoplus_{i=0}^\infty A_i$ denote a connected, $\N$-graded
$\kk$-algebra that is generated by~$A_1$.  
Suppose $A$ is Auslander-Gorenstein of finite injective dimension 
and satisfies the Cohen-Macaulay property, and that there exists a normalizing 
sequence $\{ y_1, \ldots , y_\nu\}\subset A\setminus \kk$ consisting 
of homogeneous elements such that 
GKdim$(A/\la y_1, \ldots , y_\nu \ra ) = 1$. If GKdim$(A) = n \in \N$,  then 
conditions I-IV are equivalent for any normalizing 
sequence~$F\subset A\setminus \kk^\times$ that consists of $n$~homogeneous 
elements. 
\qedsymbol
\end{thm}

Given Theorems~\ref{thm2} and \ref{thm3}, we propose the following
definition of complete intersection. 

\begin{defn}\label{compintn2}
Let $A$ be as in Theorem~\ref{thm3}.  If $F$ is a normalizing sequence in $A$
that consists of $n$ homogeneous elements of positive degree, then we call 
$\ds A/\la F \ra$ a {\em complete intersection} if the equivalent conditions 
I-IV from \S\ref{conds} hold.
\end{defn}

\noindent
Definition~\ref{compintn} is a special case of this definition, as is the 
commutative definition. By Theorem~\ref{thm2}, this definition applies to
regular GSCAs. Moreover, it applies to some other algebras as the following 
examples demonstrate.

\begin{example}
Let $A$ denote the coordinate ring of quantum $m \times \ell$ matrices
(\cite{FRT}).
By \cite[\S3]{GL}, \cite[Lemma]{Lev.Staff} and \cite[Corollary~5.10]{Lev},
$A$~satisfies the hypotheses of Theorem~\ref{thm3} (where $n = m\ell$). 
Hence, Definition~\ref{compintn2} applies to~$A$.
\end{example}

\begin{example}
Let $A$ denote the homogenization of the universal enveloping algebra of 
any Lie algebra of dimension $n-1 < \infty$. By \cite{Bruyn.VdB},
$A$~satisfies the hypotheses of Theorem~\ref{thm3}, so  
Definition~\ref{compintn2} applies to~$A$.
\end{example}

See \cite{Y} for examples of other quantized algebras that satisfy
Definition~\ref{compintn2}, such as subalgebras $U^{\pm}[w]$
(\cite{DKP,Lus}) of the quantized universal enveloping algebra
$U_q({\mathfrak{g}})$ of a simple Lie algebra $\mathfrak{g}$, where
$q\in \kk^\times$ is not a root of unity and $w$ is an element of
length~$n$ of the associated Weyl group.

\begin{ques}
Suppose $A$ and $F$ are as in \S\ref{conds} and let $I_k = \la f_1, 
\ldots , f_k\ra$ for all $k\leq t \leq n$. If $A$ is commutative, 
then, for each $k$, $\widehat{\V(I_k)}$ is a (projective) scheme, and so has 
a dimension. 
In particular, if $A$ is the polynomial ring, then $F$ is regular if and
only if 
$\dim(\widehat{\V(I_k)}) = n- k - 1$, for all $k \leq t$. However,
if $A$ is not commutative, is there an analogous statement and under what
hypotheses on $A$ could it hold?
\end{ques}



\end{document}